\documentclass[11pt, reqno]{amsart}
\usepackage{amssymb,amsmath,amscd}
\usepackage[perpage,symbol]{footmisc}
\begin{document}
\centerline{\textbf{On sums of subsets of Chen primes\footnote{This work is supported by the National Natural Science Foundation of China (Grant No. 10771135) and Innovation Program of Shanghai Municipal Education Commission.  Classification: 11P32 Keywords: Chen primes, positive relative density, W-trick, additive combinatorics}}}
\bigskip
\centerline{Zhen Cui, Hongze Li and Boqing Xue\footnote{The third author is the corresponding author.}}
\bigskip

\begin{quote}
{\footnotesize {
\begin{center}
\scshape Abstract
\end{center}
}
 \medskip

 In this paper we show that if $A$ is a subset of Chen primes with positive relative density $\alpha$, then $A+A$ must have positive upper density at least $c\alpha e^{-c^\prime\log(1/\alpha)^{2/3}(\log\log(1/\alpha))^{1/3}}$ in the natural numbers.}
\end{quote}
\bigskip

\section{ Introduction}

In 1953, K. Roth \cite{R} proved that any subset of positive integers of positive density contains non-trivial three-term arithmetic progressions. In recent years, Green \cite{G} showed that any subset of primes of relative positive density also has this property. And later Roth's theorem was extended to Chen primes in \cite{G-T}. Moreover, a celebrated theorem  was proved by Green and Tao \cite{G-T 2}, showing that the primes contain arbitrarily long arithmetic progressions.

The strategy developed by Green and Green-Tao is called ``W-trick''. The primes are embedded to a set behaving more ``pseudorandom'', meanwhile slight density-increment is gained. For various applications, one can see \cite{C-H}, \cite{CLX}, \cite{Le}, \cite{L-P}, \cite{L-P 2}.... In \cite{C-H}, it is proved by Chipeniuk and Hamel that if $A$ is a subset of primes, with positive relative density $\alpha_0$, then the set $A+A$ has positive density at least
\begin{displaymath}
C_1\alpha_0 e^{-C_2(\log(1/\alpha_0)^{2/3}(\log\log(1/\alpha_0))^{1/3})}
\end{displaymath}
in the natural numbers. This result is not far from best possible due to some examples.

Let $\mathcal{P}_c$ be the set of Chen primes, each of whom is a prime $p$ for which $p+2$ is either a prime or a product $p_1p_2$ with $p_1,p_2>p^{3/11}$, according to Chen\cite{C} and Iwaniec \cite{I}. Chen's famous theorem concludes that there are infinitely many such primes.
\medskip

\textbf{Theorem 0. } (\cite{I}) Let $n$ be a large integer. Then the number of Chen primes less than $n$ is at least $c_1n/\log^2n$, for some absolute constant $c_1>0$.

\medskip

In this paper, we extend the density result to subsets of Chen primes. For any set $S\subseteq \mathbb{N}$, denote
\begin{displaymath}
\overline{d}(S)=\limsup\limits_{n\rightarrow\infty}\frac{|S\cap [1,n]|}{|[1,n]|},\quad \overline{d}_{\mathcal{P}_c}(S)=\limsup\limits_{n\rightarrow\infty}\frac{|S\cap\mathcal{P}_c\cap [1,n]|}{|\mathcal{P}_c\cap[1,n]|}.
\end{displaymath}
\medskip

\textbf{Theorem 1.} Let $A\subseteq \mathcal{P}_c$ with positive relative density $\overline{d}_{\mathcal{P}_c}(A)=\alpha_0$. Then
\begin{displaymath}
\overline{d}(A+A)\geq C_3\alpha_0 e^{-C_4(\log(1/\alpha_0)^{2/3}(\log\log(1/\alpha_0))^{1/3})}
\end{displaymath}
for some absolute positive constant $C_3$ and $C_4$.
\medskip

Since $\overline{d}_{\mathcal{P}_c}(A)=\alpha_0$, there exist infinitely many $n$ such that $|A\cap [1,n]|/|\mathcal{P}_c\cap [1,n]|\geq \alpha_0/2$. The previous theorem will follow from a finite version, with $\alpha=\alpha_0/2$.
\medskip

\textbf{Theorem 2.} Suppose that $n$ is a sufficiently large integer. Let $A\subset \mathcal{P}_c\cap [1,n]$ with $\frac{|A|}{|\mathcal{P}_c\cap [1,n]|}\geq \alpha$. Then
\begin{displaymath}
|A+A|\geq C_5\alpha e^{-C_6\log(1/\alpha)^{2/3}(\log\log(1/\alpha))^{1/3}}n.
\end{displaymath}
for some absolute positive constant $C_5$ and $C_6$.
\medskip

We mainly follow arguments of Chipeniuk-Hamel\cite{C-H} and combine the envelop sieve function of Green-Tao\cite{G-T}.

\bigskip

\section{ Notations and Preliminary Lemmas}

For a parameter $k$, we write $f\ll_k g$ or $f=\textit{O}_k(g)$ to denote the estimate $f\leq C_k g$ for some positive constant $C_k$ depending only on $k$. For a set S, $|S|$ and $\#S$ both denote the cardinality of $S$ and the characteristic function $1_S(x)$ takes value $1$ for $x\in S$ and $0$ otherwise. The sum set $S+S^\prime:=\{s+s^\prime\;:\;s\in S,\;s^\prime \in S^\prime\}$. $c,c_0,c_1,c_2,\ldots$ are positive absolute constants. We write $\mathbb{Z}_N$ for the cyclic group $\mathbb{Z}/N\mathbb{Z}$ and $\mathbb{Z}_N^\ast$ for the multiplicative subgroup of integers modulo N. And $[1,N]$ denotes the set $\{1,2,\cdots,N\}$.

Next we introduce Fourier analysis on $\mathbb{Z}_N$. If $f:\mathbb{Z}_N\rightarrow \mathbb{C}$ is a function and $S \subseteq \mathbb{Z}_N$, we define
\begin{displaymath}
\mathbb{E}_{x\in S}f(x):=\frac{1}{|S|}\sum\limits_{x\in S}f(x).
\end{displaymath}
Let $e_N(x):=e^{2\pi i x/N}$. The Fourier transform and the Fourier inversion are the following
\begin{displaymath}
\widehat{f}(\xi)=\mathbb{E}_{x\in \mathbb{Z}_N}f(x)e_N(-\xi x), \quad\quad f(x)=\sum\limits_{\xi \in \mathbb{Z}_N} \widehat{f}(\xi)e_N(\xi x).
\end{displaymath}
The $L^q,\;L^\infty,\;l^q,\;l^\infty$-norms are defined to be
\begin{displaymath}
\|f\|_{L^q}=\left(\mathbb{E}_{x\in \mathbb{Z}_N}|f(x)|^q\right)^{1/q},\quad \|f\|_{L^\infty}=\sup\limits_{x\in \mathbb{Z}_N}|f(x)|.
\end{displaymath}
\begin{displaymath}
\|\widehat{f}\|_{l^q}=\left(\sum\limits_{\xi\in \mathbb{Z}_N}|\widehat{f}(\xi)|^q\right)^{1/q},\quad \|\widehat{f}\|_{l^\infty}=\sup\limits_{\xi\in \mathbb{Z}_N}|\widehat{f}(\xi)|.
\end{displaymath}
Plancherel's equality tells that $\|f\|_{L^2}=\|\widehat{f}\|_{l^2}$. We also write
\begin{displaymath}
f\ast g(x)=\mathbb{E}_{y\in \mathbb{Z}_N}f(x-y)g(y)
\end{displaymath}
for convolution. A basic identity for convolution is $\widehat{f\ast g}=\widehat{f}\cdot \widehat{g}$. For non-negative valued function $f$ and $g$, it obeys that
\begin{equation}
\|f\ast g\|_{L^1}=\|f\|_{L^1}\cdot \|g\|_{L^1}. \label{L1}
\end{equation}
Hausdroff-Young inequality shows that
\begin{displaymath}
\|\widehat{f}\|_{l^{q^\prime}}\leq \|f\|_{L^q},
\end{displaymath}
where $1\leq q\leq 2$ and $1/q+1/q^\prime=1$. H\"{o}lder's inequality says that
\begin{displaymath}
\|fg\|_{L^r}\leq \|f\|_{L^q}\|g\|_{L^{q^\prime}}
\end{displaymath}
whenever $0<q,q^\prime,r\leq \infty$ are such that $1/q+1/q^\prime=1/r$. And Young inequality tells that
\begin{displaymath}
\|f\ast g\|_{L^r}\leq \|f\|_{L^q}\|g\|_{L^{q^\prime}}
\end{displaymath}
for $1\leq q,q^\prime,r\leq \infty$ and $1/q+1/q^\prime=1/r+1$. (See \cite{T-V} for details.)
\medskip

Throughout this paper, $A\subseteq \mathcal{P}_c\cap[1,n]$ with $|A|/|\mathcal{P}_c\cap [1,n]|\geq \alpha$. And $p$ always denotes a prime. Let $\lambda$, $\varepsilon$, $\delta$ be small parameters and $t\gg 1$ be a very large real number to be specified later.

Write $W:=\prod\limits_{3\leq p\leq t}p$. For $b$ with $(b,W)=(b+2,W)=1$, write $F^{(b)}(x)=(Wx+b)(Wx+b+2)$, and let $R=\lfloor N^{1/20}\rfloor$. In \cite{G-T}, Green and Tao constructed an enveloping sieve function $\beta_R^{(b)}$ with the property (See \cite{G-T}Proposition 3.1)
\begin{equation}
\beta_R^{(b)}(x)\gg \mathfrak{S}_{F^{(b)}}^{-1}\log^2 R\cdot 1_{X^{(b)}_{R!}}(x) \label{beta_R}
\end{equation}
(One can check that the constant does not depend on $b$.) with
\begin{displaymath}
\mathfrak{S}_{F^{(b)}}=\prod\limits_p \frac{\gamma^{(b)}(p)}{\left(1-\frac{1}{p}\right)^2},\quad \gamma^{(b)}(p)=\frac{1}{p}|\{n\in \mathbb{Z}/p\mathbb{Z}\;:\;(p,F^{(b)}(n))=1\}|,
\end{displaymath}
and
\begin{displaymath}
X_{R!}^{(b)}=\left\{n\in \mathbb{Z}\; : \; (d,F^{(b)}(n))=1 \text{ for all }1\leq d\leq R \right\}.
\end{displaymath}

By computation (or See \cite{G-T} (6.3)),
\begin{equation}
\log^2t\ll\mathfrak{S}_{F^{(b)}}\ll\log^2 t. \label{S_F}
\end{equation}

Restrict $\beta_R^{(b)}$ to the set $[1,N]$, which we identify with $\mathbb{Z}_N$, and let $\nu^{(b)}\;:\; \mathbb{Z}_N\rightarrow\mathbb{R}^+$ be the resulting function. We have the following
\medskip

\textbf{Lemma 1.} (\cite{G-T}, Lemma 6.1) For $\xi\in \mathbb{Z}_N$, we have
\begin{displaymath}
\widehat{\nu^{(b)}}(\xi)=\delta_{\xi,0}+\textit{O}(t^{-1/2}),
\end{displaymath}
where $\delta_{\xi,0}$ is the \emph{Kronecker} delta function. Also the constant does not depend on $b$.
\medskip

\textbf{Lemma 2.} (\cite{G-T}, Prop 4.2, Lemma 3R.1) Let $q>2$ be a real number and $\{a_x\}$ be an arbitrary sequence of complex numbers. Suppose that $1\ll R\ll N^{1/10}$ and $F^{(b)},\beta_R^{(b)}$ are defined above. Then we have
\begin{displaymath}
\left(\sum\limits_{\xi\in \mathbb{Z}_N}\left|\mathbb{E}_{1\leq x\leq N}a_x\beta_R^{(b)}(x)e_N(-\xi x)\right|^q\right)^{1/q}\ll_q \left(\mathbb{E}_{1\leq x\leq N}|a_x|^2\beta_R^{(b)}(x)\right)^{1/2}
\end{displaymath}
and
\begin{displaymath}
\mathbb{E}_{1\leq x\leq N}\beta^{(b)}_R(x)\ll 1.
\end{displaymath}

\bigskip

\section{W-trick to Chen Primes}

Denote $A_n:=A\cap (\sqrt{n},n]$. And let $\Phi_W := \{b\in \mathbb{Z}_W\;:\; (b,W)=(b+2,W)=1\}$. Observe that
\begin{displaymath}
\varphi_W:= |\Phi_W|=W\prod\limits_{3\leq p \leq t}\left(1-\frac{2}{p}\right)\ll \frac{W}{\log^2 t}.
\end{displaymath}
Choose a prime $N\in (\frac{2n}{W},\frac{4n}{W}]$. For $b\in \Phi_W$, let
\begin{displaymath}
A_n^{(b)}:= \{x\in A_n\;:\;x\equiv b(\text{mod }W)\},\quad A_N^{(b)}:=\{x\leq N\;:\; Wx+b\in A_n^{(b)} \}.
\end{displaymath}
The choice of $N$ ensures that $A_n^{(b)}+A_n^{(b)}$ in $\mathbb{Z}$ can be identified with $A_N^{(b)}+A_N^{(b)}$ in $\mathbb{Z}_N$.
Noting that $A_N^{(b)}\subseteq X_{R!}^{(b)}$. Combining (\ref{beta_R}) and (\ref{S_F}), we can define $f^{(b)}:\; \mathbb{Z}_N\rightarrow \mathbb{R}^+$ by
\begin{displaymath}
f^{(b)}(x)=c\frac{\log^2N}{\log^2t}1_{A_N^{(b)}}(x),
\end{displaymath}
with $c>0$ chosen sufficiently small such that $f^{(b)}(x)\leq \nu^{(b)}(x)$. (Noting that the choice of c does not depend on $b$.)

Now follow \cite{G} and \cite{G-T}, we decompose $f^{(b)}$ into one \textsl{anti-uniform }component and one \textsl{uniform} component. Denote
\begin{displaymath}
R^{(b)}=\{\xi\in \mathbb{Z}_N\;:\; |\widehat{f^{(b)}}(\xi)|\geq \delta\},\quad B^{(b)}=\{x\in \mathbb{Z}_N: \sup\limits_{\xi\in R^{(b)}}|1-e_N(\xi x)|\leq \varepsilon\}.
\end{displaymath}
Let $\beta^{(b)}(x)=\frac{N}{|B^{(b)}|}1_{B^{(b)}}(x)$ and
\begin{displaymath}
f_1^{(b)}=f^{(b)}\ast \beta^{(b)}\ast \beta^{(b)},\quad f_2^{(b)}=f^{(b)}-f_1^{(b)}.
\end{displaymath}
\medskip

\textbf{Lemma 3.} The functions $f^{(b)},f_1^{(b)},f_2^{(b)}$ defined above have the following properties:

(\romannumeral1) $\|f_1^{(b)}\|_{L^1}=\|f^{(b)}\|_{L^1}\leq  1+\textit{O}(t^{-1/2})$.

(\romannumeral2) $\left\|\widehat{f_2^{(b)}}\right\|_{l^\infty}\ll \delta+\varepsilon\cdot(1+\textit{O}(t^{-1/2})).$

(\romannumeral3) $\left\|\widehat{f_1^{(b)}}\right\|_{l^{2+\lambda}},\;\left\|\widehat{f_2^{(b)}}\right\|_{l^{2+\lambda}}\ll\left\|\widehat{f^{(b)}}\right\|_{l^{2+\lambda}}\ll_\lambda 1$.

(\romannumeral4) $\left\|f_1^{(b)}\right\|_{L^\infty}\leq 1+\textit{O}(t^{-1/2})+\textit{O}_\lambda(\varepsilon^{-c_\lambda\delta^{-(2+\lambda)}}t^{-1/2}).$
\medskip

Proof: The proofs follow as in \cite{G-T}. We reiterate them here in order to specify $\varepsilon$ and $\delta$.

(\romannumeral 1) One can deduce that $\|\beta^{(b)}\|_{L^1}=1$. Then by (\ref{L1}) and Lemma 1, we have
\begin{displaymath}
\|f_1^{(b)}\|_{L^1}=\|f^{(b)}\|_{L^1}\leq \|\nu^{(b)}\|_{L^1}=|\widehat{\nu^{(b)}}(0)|= 1+\textit{O}(t^{-1/2}).
\end{displaymath}

(\romannumeral 2) Noting that
\begin{displaymath}
\widehat{f^{(b)}_2}{(\xi)}=\widehat{f^{(b)}}(\xi)-\widehat{f^{(b)}_1}(\xi)=\widehat{f^{(b)}}(\xi)(1-\widehat{\beta^{(b)}}(\xi)^2).
\end{displaymath}
For $\xi \in R^{(b)}$, we have
\begin{displaymath}
|1-\widehat{\beta^{(b)}}(\xi)|=|\mathbb{E}_{x\in B^{(b)}}(1-e_N(-\xi x))|\leq\mathbb{E}_{x\in B^{(b)}}|1-e_N(-\xi x)| \leq \varepsilon.
\end{displaymath}
By triangle inequality, $|1-\widehat{\beta^{(b)}}(\xi)^2|\leq 2\varepsilon$. Also we have
\begin{displaymath}
\|\widehat{f^{(b)}}\|_{l^\infty}\leq \|f^{(b)}\|_{L^1}\leq \|\nu^{(b)}\|_{L^1}=1+\textit{O}(t^{-1/2})
\end{displaymath}
by Hausdorff-Young inequality and Lemma 1. Hence
\begin{displaymath}
|\widehat{f^{(b)}_2}{(\xi)}|\ll \varepsilon\cdot (1+\textit{O}(t^{-1/2})).
\end{displaymath}

And for $\xi \notin R^{(b)}$, we have
\begin{displaymath}
|\widehat{f^{(b)}_2}{(\xi)}|\ll |\widehat{f^{(b)}}{(\xi)}|\ll\delta.
\end{displaymath}

(\romannumeral 3) For arbitrary $\lambda>0$, taking $a_x=f^{(b)}(x)/\beta_R^{(b)}(x)$ and $q=2+\lambda$, Lemma 2 tells

\begin{displaymath}
\|\widehat{f^{(b)}}\|_{l^{2+\lambda}}\ll_\lambda 1.
\end{displaymath}
Since $\|\widehat{\beta^{(b)}}\|_{l^\infty}\leq \|\beta^{(b)}\|_{L^1}= 1$, we have
\begin{displaymath}
\left|\widehat{f_1^{(b)}}(\xi)\right|=|\widehat{f^{(b)}}(\xi)||\widehat{\beta^{(b)}}(\xi)|^2\leq |\widehat{f^{(b)}}(\xi)|,
\end{displaymath}
and
\begin{displaymath}
\left|\widehat{f_2^{(b)}}(\xi)\right|\leq \left|\widehat{f^{(b)}}(\xi)\right|+\left|\widehat{f_1^{(b)}}(\xi)\right|\ll\left|\widehat{f^{(b)}}(\xi)\right|.
\end{displaymath}
Then (\romannumeral 3) follows.

(\romannumeral 4)
\begin{displaymath}
\begin{aligned}
\quad |f_1^{(b)}(x)|&=|f^{(b)}\ast \beta^{(b)}\ast \beta^{(b)}(x)|\\
&\leq |\nu^{(b)}\ast \beta^{(b)}\ast \beta^{(b)}(x)|\\
&=|\sum\limits_{\xi\in \mathbb{Z}_N}\widehat{\nu^{(b)}}(\xi)\widehat{\beta^{(b)}}(\xi)^2e_N(\xi x)|\\
&\leq |\widehat{\nu^{(b)}}(0)||\widehat{\beta^{(b)}}(0)|^2+\sup\limits_{\xi \neq 0}|\widehat{\nu^{(b)}}(\xi)|\sum\limits_{\xi}|\widehat{\beta^{(b)}}(\xi)|^2\\
&\leq (1+\textit{O}(t^{-1/2}))+\textit{O}(t^{-1/2})\cdot\frac{N}{|B^{(b)}|}.
\end{aligned}
\end{displaymath}
In the last inequality we have used Plancherel's equality to get
\begin{displaymath}
\sum\limits_{\xi}\left|\widehat{\beta^{(b)}}(\xi)\right|^2=\mathbb{E}_{x\in \mathbb{Z}_N}\left|\beta^{(b)}(x)\right|^2=\frac{N}{|B^{(b)}|}.
\end{displaymath}

By pigeonhole principle (or see \cite{Cro}, Lemma 3) we have $|B^{(b)}|\gg (\varepsilon/2\pi)^{|R^{(b)}|} N$. And from (\romannumeral3), there exists some positive integer $c_\lambda$ such that
\begin{displaymath}
c_\lambda \geq \|\widehat{f^{(b)}}\|_{l^{2+\lambda}}\geq \sum\limits_{\xi\in R}|\widehat{f}(\xi)|^{2+\lambda} \geq  |R^{(b)}|\cdot\delta^{2+\lambda},
\end{displaymath}
which implies $|R^{(b)}|\leq c_\lambda\delta^{-(2+\lambda)}$. Hence we get
\begin{displaymath}
0\leq f_1^{(b)}(x) \leq 1+\textit{O}(t^{-1/2})+\textit{O}_\lambda(\varepsilon^{-c_\lambda\delta^{-(2+\lambda)}}t^{-1/2}).
\end{displaymath}

\rightline{$\Box$}
\bigskip

Now for $b\in \Phi_W$, define
\begin{displaymath}
\delta_b:= \frac{|A_n^{(b)}|}{\frac{c_1n}{\log^2n \cdot \varphi_W}}.
\end{displaymath}
By (\romannumeral1) of Lemma 3, the definition of $f^{(b)}$ and recalling that $N\in (\frac{2n}{W},\frac{4n}{W}]$, we can get
\begin{equation}
\delta_b \ll \|f^{(b)}\|_{L^1}\ll \delta_b\ll 1, \label{fL1}
\end{equation}
when $t\gg 1$ is sufficient large. However, $\delta_b\ll 1$ is equivalent to
\begin{equation}
|A_n^{(b)}|\ll \frac{n}{\log^2n\cdot \varphi_W}. \label{Anb}
\end{equation}
\medskip

\textbf{Lemma 4.} The convolution of functions $f_1^{(b)}$ and $f_2^{(b)}$ defined above have following properties:

(\romannumeral1) $\delta_{b_1}\delta_{b_2}\ll\|f_1^{(b_1)}\ast f_1^{(b_2)}\|_{L^1}\ll \delta_{b_1}\delta_{b_2}.$

(\romannumeral2) $\|f_1^{(b_1)}\ast f_2^{(b_2)}\|_{L^2}\ll(\delta+\varepsilon\cdot(1+\textit{O}(t^{-1/2})))^{1-\lambda/2},\quad \|f_2^{(b_1)}\ast f_2^{(b_2)}\|_{L^2}\ll(\delta+\varepsilon\cdot(1+\textit{O}(t^{-1/2})))^{1-\lambda/2}.$

(\romannumeral3) $\|f_1^{(b_1)}\ast f_1^{(b_2)}\|_{L^\infty}\ll \min\{\delta_{b_1},\delta_{b_2}\}(1+\textit{O}(t^{-1/2})+\textit{O}_\lambda(\varepsilon^{-c_\lambda\delta^{-(2+\lambda)}}t^{-1/2})).$
\medskip

Proof: (\romannumeral1) By (\ref{L1}) and (\ref{fL1}), we have
\begin{displaymath}
\|f_1^{(b_1)}\ast f_1^{(b_2)}\|_{L^1}=\|f_1^{(b_1)}\|_{L^1}\|f_1^{(b_2)}\|_{L^1}=\|f^{(b_1)}\|_{L^1}\|f^{(b_2)}\|_{L^1}\gg \delta_{b_1}\delta_{b_2}.
\end{displaymath}

(\romannumeral2) Using Plancherel's equality, H\"{o}lder's inequality, and together with Lemma 3, we can obtain
\begin{displaymath}
\begin{aligned}
\|f_1^{(b_1)}\ast f_2^{(b_2)}\|_{L^2}^2&=\left\|\widehat{f_1^{(b_1)}}\widehat{ f_2^{(b_2)}}\right\|_{l^2}^2=\left\|\widehat{f_1^{(b_1)}}^2\widehat{ f_2^{(b_2)}}^2\right\|_{l^1}\\
&\leq \left\|\widehat{f_2^{(b_2)}}^{2-\lambda}\right\|_{l^\infty}\left\|\widehat{f_1^{(b_1)}}^2\widehat{f_2^{(b_2)}}^{\lambda}\right\|_{L^1} \\ &\leq \left\|\widehat{f_2^{(b_2)}}\right\|_{l^\infty}^{2-\lambda}\cdot\left\|\widehat{f_1^{(b_1)}}^2\right\|_{l^{(2+\lambda)/2}}\cdot\left\|\widehat{f_2^{(b_2)}}^{\lambda}\right\|_{l^{(2+\lambda)/\lambda}}\\
&\leq \left\|\widehat{f_2^{(b_2)}}\right\|_{l^\infty}^{2-\lambda}\cdot\left\|\widehat{f_1^{(b_1)}}\right\|_{l^{2+\lambda}}^2\cdot\left\|\widehat{f_2^{(b_2)}}\right\|_{l^{2+\lambda}}^{\lambda}\\ &\ll (\delta+\varepsilon\cdot(1+\textit{O}(t^{-1/2})))^{2-\lambda}.
\end{aligned}
\end{displaymath}

(\romannumeral3) With the Young inequality and Lemma 3, it follows that
\begin{displaymath}
\begin{aligned}
\|f_1^{(b_1)}\ast f_1^{(b_2)}\|_{L^\infty}&\leq \|f_1^{(b_1)}\|_{L^1}\cdot\|f_1^{(b_2)}\|_{L^\infty}\\
&\ll\delta_{b_1}(1+\textit{O}(t^{-1/2})+\textit{O}_\lambda(\varepsilon^{-c_\lambda\delta^{-(2+\lambda)}}t^{-1/2})).
\end{aligned}
\end{displaymath}
Similarly, we have
\begin{displaymath}
\|f_1^{(b_1)}\ast f_1^{(b_2)}\|_{L^\infty} \ll \delta_{b_2}(1+\textit{O}(t^{-1/2})+\textit{O}_\lambda(\varepsilon^{-c_\lambda\delta^{-(2+\lambda)}}t^{-1/2})).
\end{displaymath}

\rightline{$\Box$}
\bigskip

For any $b_1,b_2\in \Phi_W$, define $T^{(b_1,b_2)}=\{x\in \mathbb{Z}_N\;:\;f^{(b_1)}\ast f^{(b_2)}(x)>0\}$.
\medskip

\textbf{Proposition 1.} For any $b_1,b_2\in \Phi_W$, we have
\begin{displaymath}
|T^{(b_1,b_2)}|\gg (\delta_{b_1}+\delta_{b_2}) N,
\end{displaymath}
provided that $\varepsilon,\;\delta\ll\min\{\delta_{b_1},\delta_{b_2}\}^\frac{5}{2-\lambda}$ and $\log t\gg _\lambda \left(\frac{1}{\delta}\right)^{2+\lambda}\log\frac{1}{\varepsilon}$.
\medskip

Proof: Without loss of generality, we suppose that $\delta_{b_1}\leq \delta_{b_2}$. And suppose $\|f_1^{(b_1)}\ast f_1^{(b_2)}\|_{L^1}=c_2\delta_{b_1}\delta_{b_2}$ with $1\ll c_2\ll 1$ (Recalling Lemma 4(\romannumeral1)). Let
\begin{displaymath}
T_{1,1}=\{x\in \mathbb{Z}_N\;:\;f_1^{(b_1)}\ast f_1^{(b_2)}(x)>c_2 \delta_{b_1}\delta_{b_2}/2\},
\end{displaymath}
and
\begin{displaymath}
T_{i,j}=\{x\in \mathbb{Z}_N\;:\;|f_i^{(b_1)}\ast f_j^{(b_2)}(x)|>c_2 \delta_{b_1}\delta_{b_2}/20\}.
\end{displaymath}
for $(i,j)=(1,2),(2,1)$ or $(2,2)$. Then
\begin{displaymath}
T^{(b_1,b_2)}\supseteq T_{1,1}\cap (T_{1,2}\cup T_{2,1}\cup T_{2,2})^c.
\end{displaymath}
So
\begin{equation}
|T^{(b_1,b_2)}|\geq |T_{1,1}|-|T_{1,2}|-|T_{2,1}|-|T_{2,2}|. \label{Prop_EQ1}
\end{equation}
Combining Lemma 4, we can get
\begin{displaymath}
\begin{aligned}
c_2 \delta_{b_1}\delta_{b_2}&=\|f_1^{(b_1)}\ast f_1^{(b_2)}\|_{L^1}=\mathbb{E}_{x\in \mathbb{Z}_N}f_1^{(b_1)}\ast f_1^{(b_2)}(x)\\
&\leq \frac{1}{N}\left(|T_{1,1}|\cdot \delta_{b_1}(1+\textit{O}(t^{-1/2})+\textit{O}_\lambda(\varepsilon^{-c_\lambda\delta^{-(2+\lambda)}}t^{-1/2}))+N\cdot c_2 \delta_{b_1}\delta_{b_2}/2\right).
\end{aligned}
\end{displaymath}
We conclude that
\begin{equation}
|T_{1,1}|\gg \delta_{b_2}(1+\textit{O}(t^{-1/2})+\textit{O}_\lambda(\varepsilon^{-c_\lambda\delta^{-(2+\lambda)}}t^{-1/2}))^{-1}N. \label{Prop_EQ2}
\end{equation}
For $(i,j)=(1,2),(2,1)$ or $(2,2)$, it appears that
\begin{displaymath}
(\delta+\varepsilon\cdot(1+\textit{O}(t^{-1/2})))^{2-\lambda}\gg \|f_i^{(b_1)}\ast f_j^{(b_2)}\|_{L^2}^2\geq \frac{1}{N}|T_{i,j}|\cdot (c_2 \delta_{b_1}\delta_{b_2}/20)^2.
\end{displaymath}
\begin{equation}
|T_{i,j}|\ll(\delta+\varepsilon\cdot(1+\textit{O}(t^{-1/2})))^{2-\lambda}(\delta_{b_1}\delta_{b_2})^{-2}N. \label{Prop_EQ3}
\end{equation}
Choose $\varepsilon=\delta\ll\min\{\delta_{b_1},\delta_{b_2}\}^\frac{5}{2-\lambda}$, $t\gg 1$ and $\log t\gg c_\lambda \left(\frac{1}{\delta}\right)^{2+\lambda}\log\frac{1}{\varepsilon}$ such that
\begin{displaymath}
|T_{1,1}|>4\max\{ |T_{1,2}|,|T_{2,1}|,|T_{2,2}|\},
\end{displaymath}
then Proposition 1 follows from (\ref{Prop_EQ1}), (\ref{Prop_EQ2}) and (\ref{Prop_EQ3}).

\rightline{$\Box$}
\bigskip

\section{Sums of Subsets of Chen Primes}

Noting that $A\subseteq \mathcal{P}_c\cap [1,n]$ with $|A|\geq \alpha |\mathcal{P}_c\cap [1,n]|$ and $A_n=A\cap (\sqrt{n},n]$. By Theorem 0, we can assert that
\begin{equation}
|A_n|\geq \frac{\alpha c_1 n}{2\log^2n}. \label{Antotal}
\end{equation}
To pick out the $A_n^{(b)}$'s with `many' elements, define
\begin{displaymath}
G:= \left\{b\in \Phi_W\;:\; \delta_b \geq \frac{\alpha}{4} \right\}.
\end{displaymath}
Recall (\ref{Anb}), i.e.
\begin{displaymath}
|A_n^{(b)}|\ll \frac{n}{\log^2n\cdot \varphi_W}.
\end{displaymath}
Since
\begin{displaymath}
\frac{\alpha c_1 n}{2\log^2 n}\leq\sum\limits_{b\in \Phi_W}|A_n^{(b)}|\leq |G|\textit{O}(1)\frac{n}{\log^2 n \cdot \varphi_W}+\varphi_W\cdot \frac{\alpha}{4}\cdot \frac{c_1n}{\log^2 n \cdot \varphi_W},
\end{displaymath}
we conclude that
\begin{displaymath}
|G|\gg \alpha \varphi_W.
\end{displaymath}

Proposition 1 tells that
\begin{displaymath}
|A_n^{(b_1)}+A_n^{(b_2)}|=|A_N^{(b_1)}+A_N^{(b_2)}|\geq |T^{(b_1,b_2)}|\gg (\delta_{b_1}+\delta_{b_2}) N\gg (\delta_{b_1}+\delta_{b_2}) n/W,
\end{displaymath}
provided that we set $\varepsilon=\delta=c_3\alpha^\frac{5}{2-\lambda}$ and $\log t= c_4c_\lambda \alpha^{-\frac{5(2+\lambda)}{2-\lambda}}\log\alpha^{-1}$. \\

Denote $\Delta_x=\max\limits_{(b_1,b_2)\in G\times G\atop b_1+b_2=x}(\delta_{b_1}+\delta_{b_2})$. we have
\begin{equation}
|A_n+A_n|\geq \sum\limits_{x\in G+G} \max\limits_{(b_1+b_2)\in G\times G\atop b_1+b_2=x}|A_n^{(b_1)}+A_n^{(b_2)}|\gg \sum\limits_{x\in G+G} \Delta_x \cdot n/W. \label{A+A}
\end{equation}

By (\ref{Antotal}),
\begin{displaymath}
\sum\limits_{b\in \Phi_W}\delta_b\geq \alpha\varphi_W/2.
\end{displaymath}
Then we have
\begin{displaymath}
\sum\limits_{b\in G}\delta_b\geq \alpha\varphi_W/4.
\end{displaymath}

\begin{displaymath}
\sum\limits_{(b_1,b_2)\in G\times G}(\delta_{b_1}+\delta_{b_2})\gg\alpha\varphi_W |G|.
\end{displaymath}
For $B\subseteq \mathbb{Z}_W$ and $x\in \mathbb{Z}_W$, denote $r_B(x)=\#\{(b_1,b_2)\in G\times G\;:\; b_1+b_2=x\}$.

\medskip

\textbf{Lemma 5.} Suppose $W\in \mathbb{Z}^+$ is a sufficiently large squarefree integer. Let $\alpha>0$ and $k$ sufficiently large. And let $B\subseteq \Phi_W$ satisfy $|B|\geq \alpha \varphi_W$. Then
\begin{displaymath}
\sum\limits_{x\in\mathbb{Z}_W}r_B(x)^k\leq \frac{e^{\widetilde{C}k^3\log k}}{\alpha^2}\frac{|G|^k\varphi_W^k}{W^{k-1}}
\end{displaymath}
for some absolute constant $\widetilde{C}>0$.
\medskip

This lemma is an analog of Proposition 14 of Chipeniuk and Hamel\cite{C-H}. It can be extended to any subsets of `sieve-type' without much modification. We put the long proof in the appendix.

By Lemma 5, we conclude
\begin{displaymath}
\sum\limits_{x\in G+G} r_G(x)^k\leq \frac{e^{\textit{O}(k^3\log k)}}{\alpha^2}\frac{|G|^k\varphi_W^k}{W^{k-1}}.
\end{displaymath}

By H\"{o}lder's inequality,
\begin{displaymath}
\begin{aligned}
\alpha\varphi_W|G|&\ll\sum\limits_{x\in G+G} r_G(x)\left(\frac{1}{r_G(x)}\sum\limits_{(b_1,b_2):b_1+b_2=x}(\delta_{b_1}+\delta_{b_2})\right)\\
&\leq \sum\limits_{x\in G+G} r_G(x) \cdot \Delta_x\\
&\leq \left(\sum\limits_{x\in G+G}r_G(x)^k\right)^{1/k}\left(\sum\limits_{x\in G+G} \Delta_x^{k/(k-1)}\right)^{(k-1)/k}\\
&\ll \frac{e^{\textit{O}(k^2\log k)}}{\alpha^{2/k}}\frac{|G|\varphi_W}{W^{(k-1)/k}} \cdot\left(\sum\limits_{x\in G+G} \Delta_x^{k/(k-1)}\right)^{(k-1)/k}.
\end{aligned}
\end{displaymath}

Noting that $\Delta_x\leq 2$. Calculation shows that
\begin{equation}
\begin{aligned}
\sum\limits_{x\in G+G} \Delta_x &\gg \sum\limits_{x\in G+G} \Delta_x^{k/(k-1)} \\
&\gg \alpha^{1+\frac{3}{k-1}}e^{-\textit{O}(\frac{k^3\log k}{k-1})} W. \label{Delta_x}
\end{aligned}
\end{equation}

Combining (\ref{A+A}) and (\ref{Delta_x}), yields
\begin{displaymath}
|A_n+A_n|\gg \alpha e^{-\textit{O}(\frac{k^3\log k}{k-1})+\frac{3\log \alpha}{k-1}} n.
\end{displaymath}
If $\alpha$ is small enough, it can be deduced that
\begin{displaymath}
|A+A|\geq |A_n+A_n|\gg \alpha e^{-\textit{O}((\log \alpha^{-1})^{2/3}(\log\log\alpha^{-1})^{1/3})} n
\end{displaymath}
by taking $k=\lfloor(\log\alpha^{-1}/\log\log\alpha^{-1})^{1/3}\rfloor$. For $\alpha$ is not small, Theorem 2 can also follow by partition $B$ into the union of smaller subsets such that the above argument can be applied. (See \cite{C-H})

\rightline{$\Box$}

\section{Appendix: Addition in $\Phi_W$}

Proof of Lemma 5: Let
\begin{displaymath}
\begin{aligned}
R(x):= &\#\{(b,r)\in B\times \Phi_W:b+r=x \}\\
=&\#\{b\in B: (b-x)(b-x+2)\not\equiv 0 (\text{mod }p) \text{ for all }p|W\}.
\end{aligned}
\end{displaymath}
Put
\begin{displaymath}
\begin{aligned}
X_d:=&\{x\in [0,W-1]:(x,W)=d\}\\
=&\{x\in [0,W-1]:x=dl \text{ for some }l\in [0,W/d-1] \text{ with } (l,W/d)=1\}.
\end{aligned}
\end{displaymath}
We have
\begin{displaymath}
\begin{aligned}
S=&\sum\limits_{x\in \mathbb{Z}_W} r_B(x)^k\\
\leq &\sum\limits_{x\in \mathbb{Z}_W} R(x)^k\\
=&\sum\limits_{d|W}\sum\limits_{x\in X_d} R(x)^k\\
=&\sum\limits_{d|W}\sum\limits_{x\in X_d}\#\{b\in B: (b-x)(b-x+2)\not\equiv 0 (\text{mod }p) \text{ for all }p|W/d\}^k\\
\leq &\sum\limits_{d|W}\sum\limits_{b_1,\ldots,b_k\in B}\#\{x\in X_d: (b_i-x)(b_i-x+2)\not\equiv 0 (\text{mod }p) \text{ for all }p|W/d\text{ and }1\leq i\leq k\}\\
\leq &\sum\limits_{d|W}\sum\limits_{b_1,\ldots,b_k\in B}\#\{l\in [0,W/d-1]: (l,W/d)=1, \\
&\quad\quad\quad (b_id^{-1}-l)(b_id^{-1}-l+2d^{-1})\not\equiv 0 (\text{mod }p) \text{ for all }p|W/d\text{ and }1\leq i\leq k\}\\
\leq &\sum\limits_{d|W}\sum\limits_{b_1,\ldots,b_k\in B}\#\{l\in [0,W/d-1]: l\not\equiv 0 \text{ for all }p|W/d, \\
&\quad\quad\quad (b_id^{-1}-l)(b_id^{-1}-l+2d^{-1})\not\equiv 0 (\text{mod }p) \text{ for all }p|W/d\text{ and }1\leq i\leq k\}\\
=&\sum\limits_{d|W}\sum\limits_{b_1,\ldots,b_k\in B}\prod_{p|W/d} (p-r_p(b_1,\ldots,b_k)-1)\\
=&\sum\limits_{d|W}\sum\limits_{b_1,\ldots,b_k\in B}\frac{W}{d}\prod_{p|W/d}\left(1-\frac{r_p(b_1,\ldots,b_k)+1}{p}\right),
\end{aligned}
\end{displaymath}
where
\begin{displaymath}
r_p(b_1,\ldots,b_k)=\#\{s\in [0,p-1]: (b_id^{-1}-s)(b_id^{-1}-s+2d^{-1})\equiv 0 (\text{mod }p) \text{ for some }1\leq i\leq k\}.
\end{displaymath}

Now we fix a $d|W$.

\textbf{Claim.} For $(b_1,\ldots,b_k)\in B^k$, let
\begin{displaymath}
f(b_1,\ldots,b_k)=\sum\limits_{p|W\atop r_p(b_1,\ldots,b_k)\leq 2k-1}\frac{1}{p}.
\end{displaymath}
and
\begin{displaymath}
K=\{(b_1,\ldots,b_k)\in B^k: f(b_1,\ldots,b_k)\geq \beta\}.
\end{displaymath}
Then there exists an absolute constant $c>0$ such that
\begin{displaymath}
|K|\leq k^22^{-\exp(\beta/ck^2)}|B|^{k-2}\varphi_W^2.
\end{displaymath}
holds uniformly for every $\beta>0$.

Proof: Given $(b_1,\ldots,b_k)\in K$, we have
\begin{displaymath}
\begin{aligned}
\beta \leq &\sum\limits_{p|W\atop r_p(b_1,\ldots,b_k)\leq 2k-1}\frac{1}{p}\\
\leq &\sum\limits_{p|W\atop p|(b_i-b_j)(b_i-b_j-2) \text{ for some }i\neq j}\frac{1}{p}\\
\leq &\sum\limits_{\{i,j\}\subseteq \{1,\ldots,k\}\atop i\neq j}\sum\limits_{p|W\atop p|(b_i-b_j)(b_i-b_j-2)}\frac{1}{p}.
\end{aligned}
\end{displaymath}
By the pigeon hole principle, there exists some pair $\{i,j\}$ with $i\neq j$ such that
\begin{displaymath}
\frac{\beta}{
k(k-1)/2}\leq \sum\limits_{p|W\atop p|(b_i-b_j)(b_i-b_j-2)}\frac{1}{p}.
\end{displaymath}
Since each $(b_1,\ldots,b_k)\in K$ contribute at least one such pair $\{i,j\}$ and a given $\{b_i,b_j\}$ comes from at most $k^2|B|^{k-2}$ $k$-tuples, hence we can conclude
\begin{equation}
\frac{|K|}{k^2|B|^{k-2}}2^l\beta^l
k^{-l}(k-1)^{-l}\leq \sum\limits_{\{b,c\}\subseteq B \atop b\neq c}\left(\sum\limits_{p|W\atop p|(b-c)(b-c-2)}\frac{1}{p}\right)^l. \label{eqeq}
\end{equation}
Furthermore,  we have
\begin{displaymath}
\begin{aligned}
&\sum\limits_{\{b,c\}\subseteq B \atop b\neq c}\left(\sum\limits_{p|W\atop p|(b-c)(b-c-2)}\frac{1}{p}\right)^l\\
\leq &\sum\limits_{(b,c)\in B^2}\left(\sum\limits_{p|W\atop p|(b-c)(b-c-2)}\frac{1}{p}\right)^l\\
=&\sum\limits_{p_1,\ldots,p_l|W}\frac{1}{p_1\ldots p_l}\sum\limits_{(b,c)\in B^2\atop (b-c)(b-c-2)\equiv 0 (\text{mod }lcm[p_1,\ldots,p_l])}1\\
\end{aligned}
\end{displaymath}
Write $p_0=\max\limits_{1\leq i\leq l}p_i$ for fixed $p_1,\ldots,p_l$, one can deduce that
\begin{displaymath}
\begin{aligned}
&\sum\limits_{(b,c)\in B^2\atop (b-c)(b-c-2)\equiv 0 (\text{mod }lcm[p_1,\ldots,p_l])}1\leq \sum\limits_{b\in \Phi_W}\sum\limits_{c\in \Phi_W \atop (b-c)(b-c-2)\equiv 0 (\text{mod }p_0)}1\\
&\quad\quad\leq 2\sum\limits_{b\in \Phi_W}\max\limits_{a\in \mathbb{Z}_W}\left\{\sum\limits_{c\in \varphi_W \atop c\equiv a (\text{mod }p_0)}1\right\}\leq 2\varphi_W\max\limits_{a\in \mathbb{Z}_W}\left\{\sum\limits_{c\equiv a (\text{mod }p_0)\atop c\not\equiv 0,-2(\text{mod }p) \text{ for all }p|W/p_0}1\right\}\\
&\quad\quad\leq 2\varphi_W\varphi_{W/P_0}\leq \frac{2\varphi_W^2}{p_0-2}\leq \frac{2\varphi_W^2}{\prod\limits_{1\leq i\leq l}(p_i-2)^{1/l}}.
\end{aligned}
\end{displaymath}
Then
\begin{displaymath}
\begin{aligned}
&\sum\limits_{\{b,c\}\subseteq B \atop b\neq c}\left(\sum\limits_{p|W\atop p|(b-c)(b-c-2)}
\frac{1}{p}\right)^l\\
\leq &2\varphi_W^2\sum\limits_{p_1,\ldots,p_l|W}\frac{1}{\prod\limits_{1\leq i\leq l}p_i(p_i-2)^{1/l}}\\
=&2\varphi_W^2\left(\sum\limits_{p|W}\frac{1}{p(p-2)^{1/l}}\right)^l\\
\leq &2\varphi_W^2\left(\sum\limits_{p\leq l^l}\frac{1}{p}+\sum\limits_{n\geq l^l}\frac{1}{n(n-2)^{1/l}}\right)^l\\
\leq &  \varphi_W^2 (c\log l)^l
\end{aligned}
\end{displaymath}
for some absolute constant $c>0$. Combining (\ref{eqeq}) and the above formula, gives
\begin{displaymath}
|K|\leq 2^{-l}\beta^{-l}k^{l+2}(k-1)^lc^l(\log l)^l |B|^{k-2}\varphi_W^2.
\end{displaymath}
The Claim follows by taking $l=\exp(\beta/ck^2)$.

\rightline{$\Box$}

\medskip

Writing
\begin{displaymath}
W_1=\prod\limits_{p|W\atop p\leq 5k} p, W_2=\prod\limits_{p|W\atop p>5k}p.
\end{displaymath}
\begin{displaymath}
d_1=\prod\limits_{p|d\atop p\leq 5k} p, d_2=\prod\limits_{p|d\atop p>5k}p.
\end{displaymath}

The estimate below will be useful later.
\begin{displaymath}
\begin{aligned}
&\log \left(\prod\limits_{p|W_2/d_2\atop r_p(b_1,\ldots,b_k)\leq 2k-1}\left(1-\frac{2k+1}{p}\right)^{-1}\right)\\
=&-\sum\limits_{p|W_2/d_2\atop r_p(b_1,\ldots,b_k)\leq 2k-1}\log\left(1-\frac{2k+1}{p}\right)\\
=&\sum\limits_{p|W_2/d_2\atop r_p(b_1,\ldots,b_k)\leq 2k-1}\frac{2k+1}{p}\sum\limits_{t=0}^\infty\frac{1}{t+1}\left(\frac{2k+1}{p}\right)^t\\
\leq &2(2k+1)\sum\limits_{p|W_2/d_2\atop r_p(b_1,\ldots,b_k)\leq 2k-1}\frac{1}{p}\\
\leq &2(2k+1)f(b_1,\ldots,b_k).
\end{aligned}
\end{displaymath}
The last step is resulted from the fact that $\frac{1}{t+1}\left(\frac{2k+1}{p}\right)^t\leq \frac{1}{2^t}$ for $t\geq 1$ and $p>5k$. Noting that $r_p(b_1,\ldots,b_k)\geq 1$. We continue to estimate $S$.
\begin{displaymath}
\begin{aligned}
S\leq &\sum\limits_{d|W}\sum\limits_{b_1,\ldots,b_k\in B}\frac{W}{d}\prod_{p|W/d}\left(1-\frac{r_p(b_1,\ldots,b_k)+1}{p}\right)\\
\leq &\sum\limits_{d|W}\sum\limits_{b_1,\ldots,b_k\in B}\frac{W_1}{d_1}\prod_{p|W_1/d_1}\left(1-\frac{2}{p}\right)\frac{W_2}{d_2}\prod_{p|W_2/d_2}\left(1-\frac{r_p(b_1,\ldots,b_k)+1}{p}\right)\\
=&\sum\limits_{d|W} \varphi_{W_1/d_1}\sum\limits_{b_1,\ldots,b_k\in B}\frac{W_2}{d_2}\prod_{p|W_2/d_2}\left(1-\frac{r_p(b_1,\ldots,b_k)+1}{p}\right)\\
\leq &\sum\limits_{d|W}\varphi_{W_1/d_1} \sum\limits_{j=-\infty}^\infty \sum\limits_{b_1,\ldots,b_k\in B \atop f(b_1,\ldots,b_k)\in [2^j,2^{j+1})}\frac{W_2}{d_2}\prod_{p|W_2/d_2}\left(1-\frac{2k+1}{p}\right)\prod\limits_{p|W_2/d_2\atop r_p(b_1,\ldots,b_k)\leq 2k-1}\left(1-\frac{2k+1}{p}\right)^{-1}\\
\leq &\sum\limits_{d|W}\varphi_{W_1/d_1} \frac{W_2}{d_2}\prod_{p|W_2/d_2}\left(1-\frac{2k+1}{p}\right) \sum\limits_{j=-\infty}^\infty \sum\limits_{b_1,\ldots,b_k\in B \atop f(b_1,\ldots,b_k)\in [2^j,2^{j+1})}e^{(2k+1)2^{j+2}}\\
\leq &\sum\limits_{d|W}k^2|B|^{k-2}\varphi_W^2\varphi_{W_1/d_1} \frac{W_2}{d_2}\prod_{p|W_2/d_2}\left(1-\frac{2k+1}{p}\right) \sum\limits_{j=0}^\infty 2^{-\exp(2^j/ck^2)}e^{(2k+1)2^{j+2}}\\
=&C_kk^2|B|^{k-2}\varphi_W^2\sum\limits_{d|W}\varphi_{W_1/d_1} \frac{W_2}{d_2}\prod_{p|W_2/d_2}\left(1-\frac{2k+1}{p}\right),
\end{aligned}
\end{displaymath}
where
\begin{displaymath}
\begin{aligned}
C_k=&\sum\limits_{j=0}^\infty 2^{-\exp({2^j}/ck^2)}e^{(2k+1)2^{j+2}}\\
=&\sum\limits_{j=0}^\infty e^{(2k+1)2^{j+2}-\exp({2^j}/ck^2)\log 2}\\
\leq &\sum\limits_{j=0}^\infty e^{(2k+1)2^{j+2}-\log 2\left(\frac{{2^j}}{ck^2}+\frac{2^{2j}}{2c^2k^4}\right)}.
\end{aligned}
\end{displaymath}
The range of summation in $j$ over $(-\infty,0)$ can be omitted since $f(b_1,\ldots,b_k)\geq \sum\limits_{p\leq 2k}\gg \log\log k>1$ for sufficiently large $k$. Now for $j\geq j_1:=\frac{\log\left( 8c^2k^4(2k+1)/\log 2\right)}{\log 2}$, we have $\frac{2^{2j}\log 2}{2c^2k^4}\geq (2k+1)2^{j+2}$ and $\frac{2^j}{ck^2}\geq 2^{j/2}\geq j$. Then
\begin{displaymath}
\begin{aligned}
\sum\limits_{j\geq j_1} e^{(2k+1)2^{j+2}-\log 2\left(\frac{{2^j}}{ck^2}+\frac{{2^j}^2}{2c^2k^4}\right)}
\leq\sum\limits_{j\geq j_1}e^{-\log 2\frac{{2^j}}{ck^2}}\leq \sum\limits_{j\geq j_1}2^{-j}\leq 1.
\end{aligned}
\end{displaymath}
For $j\leq j_1$, the exponent $(2k+1)2^{j+2}-\exp({2^j}/ck^2)\log 2$ is maximized when
\begin{displaymath}
2^j=ck^2\log(4(2k+1)ck^2/\log 2)
\end{displaymath}
and
\begin{displaymath}
\begin{aligned}
&\sum\limits_{j=0}^{j_1} e^{(2k+1)2^{j+2}-\exp({2^j}/ck^2)\log 2} \\
& \quad\quad \leq  \left(\frac{\log\left( 8c^2k^4(2k+1)/\log 2\right)}{\log 2}\right)\cdot e^{2(2k+1)ck^2\left(2\log(2(2k+1)ck^2/\log 2)\right)}
\end{aligned}
\end{displaymath}
Now we get
\begin{displaymath}
C_k\leq e^{c_5k^3\log k}
\end{displaymath}
for some absolute constant $c_5>0$.

Now we turn back to $S$. The number of divisors $d$ of $W$ which gives the same $d_2$ is smaller than $\sum\limits_{t=0}^{5k} \binom{5k}{t}=2^{5k}$. And recall $|B|=\alpha \varphi_W$. Hence
\begin{displaymath}
\begin{aligned}
S\leq &C_kk^2|B|^{k-2}\varphi_W^2\sum\limits_{d|W}\varphi_{W_1/d_1} \frac{W_2}{d_2}\prod_{p|W_2/d_2}\left(1-\frac{2k+1}{p}\right)\\
\leq &C_kk^22^{5k}|B|^{k-2}\varphi_W^2W\sum\limits_{d_2|W_2} \frac{1}{d_2}\prod_{p|W_2/d_2}\left(1-\frac{2k+1}{p}\right)\\
\leq &C_k^2\alpha^{-2}|B|^kW\prod\limits_{p|W_2}\left(1-\frac{2k+1}{p}\right)\cdot\sum\limits_{d_2|W_2}\frac{1}{d_2}\prod\limits_{p|d_2}\left(1-\frac{2k+1}{p}\right)^{-1}\\
\leq &C_k^2\alpha^{-2}|B|^kW\prod\limits_{p|W_2}\left(1-\frac{2k+1}{p}\right)\left(1+\frac{1}{p-2k-1}\right)\\
\leq &C_k^2\alpha^{-2}|B|^kW\prod\limits_{p|W_2}\left(1-\frac{2}{p}\right)^{k}\\
=&C_k^2\alpha^{-2}|B|^k\varphi_W^kW^{-k+1}\prod\limits_{p|W_1}\left(1-\frac{2}{p}\right)^{-k}\\
\leq &C_k^3\alpha^{-2}|B|^k\varphi_W^kW^{-k+1}\\
\leq &\frac{e^{c_6k^3\log k}|B|^k\varphi_W^k}{\alpha^2W^{k-1}}.
\end{aligned}
\end{displaymath}

\rightline{$\Box$}

\end{document}